%% file: m4-12.tex
\documentclass{gtart}

\input gtmonout

\volumenumber{4}
\volumename{Invariants of knots and 3-manifolds (Kyoto 2001)}
\volumeyear{2002}
\papernumber{12}
\pagenumbers{183}{199}
\received{19 December 2001}
\revised{15 February 2002}
\accepted{22 July 2002}
\published{21 September 2002}

\usepackage{amsmath,amssymb,pstricks}

\newtheorem{thm}{Theorem}[section]

\newtheorem{proposition}[thm]{Proposition}

\theoremstyle{definition}

\newtheorem{definition}[thm]{Definition}

\newcommand{\CC}{{\mathbb C}}

\newcommand{\NN}{{\mathbb N}}
\newcommand{\QQ}{{\mathbb Q}}
\newcommand{\RR}{{\mathbb R}}

\newcommand{\ZZ}{{\mathbb Z}}
\newcommand{\CA}{{\cal A}}

\newcommand{\bp}{\begin{proof}}
\newcommand{\eop}{\end{proof}}

\newcommand{\tata}{\begin{pspicture}[0.2](0,0)(.5,.4)
\pscircle(0.25,0.2){.25}
\psline[linestyle=dashed,dash=2pt 1pt]{*-*}(0.05,.2)(.45,.2)
\end{pspicture}}

\newcommand{\fcpo}{\begin{pspicture}[.2](0,0)(.6,.45)
\psline{->}(0.05,0)(.45,.4)
\psline{-}(.45,0)(.3,.15)
\psline{->}(.2,.25)(0.05,.4)
\end{pspicture}}

\newcommand{\fcmo}{\begin{pspicture}[.2](0,0)(.6,.45)
\psline{->}(0.45,0)(.05,.4)
\psline{-}(.05,0)(.2,.15)
\psline{->}(.3,.25)(0.45,.4)
\end{pspicture}
}

\newcommand{\trid}{\begin{pspicture}[.2](0,0)(.3,.4)
\psline{->}(0.05,0)(.05,.4)
\psline{<-}(.25,0)(.25,.4)
\end{pspicture}
}

\newcommand{\smaxo}{\begin{pspicture}[.2](0,0)(.5,.4)
\pscurve{->}(.05,0)(0.25,0.4)(.45,0)
\end{pspicture}
}

\newcommand{\smaxno}{\begin{pspicture}[.2](0,0)(.5,.6)
\pscurve{->}(.3,.25)(0.45,.4)(.25,.6)(.05,.4)(.25,.2)(0.45,0)
\psline{-}(.05,0)(.2,.15) 
\end{pspicture} }

\begin{document}
\title{On configuration space integrals for links}
\author{Christine Lescop}
\address{CNRS, Institut Fourier, B.P.74, 38402 Saint-Martin-d'H\`eres cedex,
France}
\asciiaddress{CNRS, Institut Fourier, B.P.74, 38402 Saint-Martin-d'Heres cedex,
France}
\email{lescop@ujf-grenoble.fr}
\url{http://www-fourier.ujf-grenoble.fr/\char'176lescop/}

\begin{abstract}
We give an introductory survey on the universal Vassiliev invariant
called {\em the perturbative series expansion of the Chern--Simons
theory of links in euclidean space\/}, and on its relation with the
Kontsevich integral.  We also prove an original geometric property of
the {\em anomaly} of Bott, Taubes, Altschuler, Freidel and
D.~Thurston, that allowed Poirier to prove that the Chern--Simons
series and the Kontsevich integral coincide up to degree 6.
\end{abstract}

\asciiabstract{ 
We give an introductory survey on the universal Vassiliev invariant
called the perturbative series expansion of the Chern-Simons theory of
links in euclidean space, and on its relation with the Kontsevich
integral.  We also prove an original geometric property of the anomaly
of Bott, Taubes, Altschuler, Freidel and D. Thurston, that allowed
Poirier to prove that the Chern-Simons series and the Kontsevich
integral coincide up to degree 6.}

\primaryclass{57M27}
\secondaryclass{57M25, 17B37, 81T18}

\keywords{Kontsevich Integral, Chern--Simons theory, Vassiliev
invariants, links, knots, tangles, configuration spaces, quantum
invariants, Jacobi diagrams}
\asciikeywords{Kontsevich Integral, Chern-Simons theory, Vassiliev
invariants, links, knots, tangles, configuration spaces, quantum
invariants, Jacobi diagrams}
\maketitle

\section{Introduction}
There are essentially two universal Vassiliev invariants of links,
the {\em Kontsevich integral,\/} and the {\em perturbative series expansion of the Chern--Simons theory\/} studied by Guadagnini Martellini and Mintchev~\cite{gua}, Bar-Natan~\cite{bn'}, Axelrod and Singer~\cite{as1,as2}, Kontsevich~\cite{Kon}, Polyak and Viro~\cite{pv}, Bott and Taubes~\cite{bt}, Altschuler and Freidel~\cite{af}, D. Thurston~\cite{th}, Yang, Poirier~\cite{Po}\dots\ 
The question that was raised by Kontsevich in \cite{Kon} whether the two invariants
coincide or not is still open though it is predicted by the Chern--Simons gauge theory that the two invariants should coincide.
In \cite{Po}, Sylvain Poirier reduced this question to the computation of the {\em anomaly} of Bott, Taubes, Altschuler, Freidel and D.~Thurston, which is an element
$\alpha$ of the space of Jacobi diagrams $\CA(S^1)$.

In this article, we shall begin with an elementary introduction to  {\em the perturbative series expansion of the Chern--Simons theory of links in euclidean space\/} defined by means of configuration space integrals in a natural and beautiful way. We shall call this series the Chern--Simons series, and we shall denote it by $Z_{CS}$. Its physical interpretation will not be treated here and we refer the reader to the survey~\cite{lab} of Labastida for the interpretation of 
$Z_{CS}$ in the context of the Chern--Simons gauge theory.

Then we shall give a short survey of the substantial Poirier work \cite{Po} that allowed the author \cite{Les} to define the isomorphism of $\CA$ which 
transforms the Kontsevich integral into the Poirier limit
of the Chern--Simons invariant of framed links, as an explicit function of $\alpha$, 
and to prove the algebraic property of the anomaly~: 
{\em The anomaly has two legs.\/} 
We end up the article by proving an additional original geometric property (Proposition~\ref{mainprop}) of the anomaly that allowed Poirier
to compute the anomaly up to degree 6.

I thank St\'ephane Guillermou, Lucien Guillou, Sylvain Poirier and especially the referee for useful comments 
on these notes. 
I also thank Tomotada Ohtsuki and Hitoshi Murakami for organizing the very interesting Kyoto conference of September 2001 at the RIMS, and for inviting me to participate.

\section{Introduction to configuration space integrals: The Gauss integrals}
\label{secgauss}

In 1833, Carl Friedrich Gauss defined the first example of a {\em configuration space integral\/} for an oriented two-component link. Let us formulate his definition in a modern 
language. Consider a smooth ($C^{\infty}$) embedding 
$$L: S^1_1 \sqcup S^1_2 \hookrightarrow \RR^3$$
of the disjoint union of two circles $S^1=\{z \in \CC \;\;\mbox{s.t.}\;\; |z|=1 \}$ 
into $\RR^3$. With an element $(z_1,z_2)$ of $S^1_1 \times S^1_2$ that will be called a {\em configuration\/}, we associate the oriented direction
$$\Psi((z_1,z_2))=\frac{1}{\parallel\overrightarrow{L(z_1)L(z_2)}\parallel}\overrightarrow{L(z_1)L(z_2)} \in S^2$$ of the vector $\overrightarrow{L(z_1)L(z_2)}$.
Thus, we have associated a map $$\Psi: S^1_1 \times S^1_2 \longrightarrow S^2$$ from a compact oriented 2-manifold to another one
with our embedding.
This map has an integral degree $\mbox{deg}(\Psi)$ that can be defined in several equivalent ways.
For example, it is the {\em differential degree\/} $\mbox{deg}(\Psi,y)$ of any regular value $y$ of $\Psi$, that is the sum of the $\pm 1$ signs of the Jacobians of $\Psi$ at the points of the preimage of $y$~\cite[\S 5]{mil}. Thus, $\mbox{deg}(\Psi)$ can easily be computed from a regular diagram of our two-component link as the differential degree of a unit vector $\overrightarrow{v}$ pointing to the reader or as the differential degree of $(-\overrightarrow{v})$.
$$\mbox{deg}(\Psi)= \mbox{deg}(\Psi,\overrightarrow{v})=  \sharp \begin{pspicture}[.2](0,0)(.6,.45)
\psline{->}(0.45,0)(.05,.4)
\rput(.25,0){\tiny 2}
\rput(.55,0){\tiny 1}
\psline[border=1pt]{->}(.05,0)(0.45,.4)
\end{pspicture} - \sharp \begin{pspicture}[.2](0,0)(.6,.45)
\rput(.25,0){\tiny 1}
\rput(.55,0){\tiny 2}
\psline{->}(.05,0)(0.45,.4)
\psline[border=1pt]{->}(0.45,0)(.05,.4)
\end{pspicture}= \mbox{deg}(\Psi,-\overrightarrow{v})=\sharp \begin{pspicture}[.2](0,0)(.6,.45)
\psline{->}(0.45,0)(.05,.4)
\rput(.25,0){\tiny 1}
\rput(.55,0){\tiny 2}
\psline[border=1pt]{->}(.05,0)(0.45,.4)
\end{pspicture} - \sharp \begin{pspicture}[.2](0,0)(.6,.45)
\rput(.25,0){\tiny 2}
\rput(.55,0){\tiny 1}
\psline{->}(.05,0)(0.45,.4)
\psline[border=1pt]{->}(0.45,0)(.05,.4)
\end{pspicture} $$
It can also be defined as the following {\em configuration space integral\/}
$$\mbox{deg}(\Psi)=\int_{S^1 \times S^1}\Psi^{\ast}(\omega)$$
where $\omega$ is the homogeneous volume form on $S^2$ such that $\int_{S^2}{\omega}=1$.
Of course, this integral degree is an isotopy invariant of $L$, and the reader has recognized that $\mbox{deg}(\Psi)$
is nothing but the {\em linking number\/} of the two components of $L$.

We can again follow Gauss and associate the following similar {\em Gauss integral\/}
$I(K;\theta)$
to a $C^{\infty}$ embedding $K: S^1 \hookrightarrow \RR^3$. (The meaning of $\theta$ will be specified later.)
Here, we consider the {\em configuration space\/} 
$C(K;\theta)=S^1 \times ]0,2\pi [$, and the map 
$$\Psi : C(K;\theta) \longrightarrow S^2$$ 
that maps $(z_1,\eta)$ to the oriented direction of $\overrightarrow{K(z_1)K(z_1e^{i\eta})}$, and we set 
$$I(K;\theta )=\int_{C(K;\theta)}\Psi^{\ast}(\omega).$$
Let us compute $I(K;\theta)$ in some cases. First notice that $\Psi$ may be extended to the closed annulus 
$$\overline{C}(K;\theta)=S^1 \times [0,2\pi]$$ by the tangent map $K^{\prime}$ of $K$ along $S^1 \times \{0\}$ and by $(-K^{\prime})$  along $S^1 \times \{2\pi\}$. Then by definition, $I(K;\theta)$ is the algebraic area (the integral of the differential degree with respect to the measure associated with $\omega$) of the image of the annulus in $S^2$.
Now, assume that $K$ is contained in a horizontal plane except in a neighborhood of crossings where it entirely lies in vertical planes. Such a knot embedding
will be called {\em almost horizontal.\/}
In that case, the image of the annulus boundary has the shape of the following bold line in $S^2$.
$$
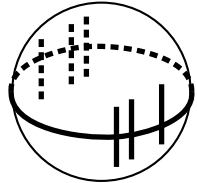
\begin{pspicture}[.4](0,0)(2.8,2.8)
\psset{unit=2cm}
\pscircle(.7,.7){.6}
\psecurve[linewidth=2pt,linestyle=dashed,dash=3pt 2pt]{-}(.7,.5)(.1,.7)(.6,1)(1.3,.7)(.7,.5)
\psline[linewidth=2pt,linestyle=dashed,dash=3pt 2pt]{-}(.6,.8)(.6,1.2)
\psline[linewidth=2pt,linestyle=dashed,dash=3pt 2pt]{-}(.5,.75)(.5,1.15)
\psline[linewidth=2pt,linestyle=dashed,dash=3pt 2pt]{-}(.3,.65)(.3,1.05)
\psecurve[linewidth=2pt]{-}(.7,.9)(.1,.7)(.8,.4)(1.3,.7)(.7,.9)
\psline[linewidth=2pt]{-}(.8,.2)(.8,.6)
\psline[linewidth=2pt]{-}(.9,.25)(.9,.65)
\psline[linewidth=2pt]{-}(1.1,.35)(1.1,.75)
\end{pspicture}$$
In particular, for each hemisphere, the differential degree of a regular value of $\Psi$ does not depend on the choice of the regular value in the hemisphere.
Assume that the orthogonal projection onto the horizontal plane is regular.
Then $I(K;\theta)$ is the average of the differential degrees of the North Pole
and the South Pole, and it can be computed from the horizontal projection as
$$I(K;\theta)=\sharp \fcpo - \sharp \fcmo.$$
This number, that is called the {\em writhe\/} of the projection, can be changed
without changing the isotopy class of the knot by local modifications where 
$\begin{pspicture}[.4](-.1,.1)(.9,.45)
\psline(0,.4)(.8,.4)\end{pspicture}$ becomes 
$\begin{pspicture}[.4](-.1,.1)(.9,.45)
\psecurve(1.8,.4)(.8,.4)(.55,.4)(.25,.25)(.4,.1)(.55,.25)
\psecurve[border=1.5pt](.25,.25)(.4,.1)(.55,.25)(.25,.4)(0,.4)(-1,.4)\end{pspicture}$ or $\begin{pspicture}[.4](-.1,.1)(.9,.45)
\psecurve(-1,.4)(0,.4)(.25,.4)(.55,.25)(.4,.1)(.25,.25)
\psecurve[border=1.5pt](.55,.25)(.4,.1)(.25,.25)(.55,.4)(.8,.4)(1.8,.4)\end{pspicture}$. In particular, $I(K;\theta)$ can reach any integral value on a given isotopy class of knots, and since it varies continuously on such a class, it can reach any real value on any given isotopy class of knots.
Thus, this Gauss integral is NOT an isotopy invariant. 

However, we can follow Guadagnini, Martellini, Mintchev~\cite{gua} and Bar-Natan~\cite{bn'}
and associate 
configuration space integrals to any embedding $L$ of an oriented one-manifold $M$ and to any Jacobi diagram $\Gamma$ on $M$. 
Let us first recall what a Jacobi diagram on a one-manifold is.

\section{Definitions of the spaces of Jacobi diagrams}
\label{secdia}

\begin{definition}
\label{defdia}
Let $M$ be an oriented one-manifold.
A {\em  Jacobi diagram\/} $\Gamma$ with support $M$ is a finite uni-trivalent graph $\Gamma$ without simple loop like $\begin{pspicture}[.2](0,0)(.6,.4)
\psline[linestyle=dashed,dash=2pt 1pt]{-*}(0.05,.2)(.25,.2)
\psline[linestyle=dashed,dash=2pt 1pt]{-}(.25,.2)(.4,.05)(.55,.2)(.4,.35)(.25,.2)
\end{pspicture}$ such that
every connected component of $\Gamma$ has at least one univalent vertex, 
equipped with:
\begin{enumerate}
\item an isotopy class of injections $i$ of the set $U$ of univalent vertices of $\Gamma$ also called {\em legs} of $\Gamma$ into the interior of $M$,
\item an {\em orientation\/} of every trivalent vertex, that is a cyclic order
on the set of the three half-edges which meet at this vertex.
\end{enumerate}
\end{definition}

Such a diagram $\Gamma$ is represented by a planar immersion of $\Gamma \cup M$ where the univalent vertices of $\Gamma$ are identified with their images under $i$, the one-manifold $M$ is represented by solid lines, whereas
the diagram $\Gamma$ is dashed. The vertices are represented by big points. The orientation of a 
vertex is represented by the counterclockwise order of the three dashed half-edges
that meet at that vertex.
Here is an example of a diagram $\Gamma$ on the disjoint 
union $M=S^1 \sqcup S^1$ of two circles:
$$
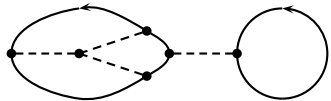
\begin{pspicture}[.2](0,0)(4.5,1.5)
\psset{unit=1.5cm}
\psecurve{->}(1.3,.7)(.7,.9)(.1,.5)(.7,.1)(1.3,.3)(1.5,.5)(1.3,.7)(.7,.9)(.1,.5)
\psarc{->}(2.5,.5){.4}{-90}{90}
\psarc{-}(2.5,.5){.4}{90}{-90}
\psline[linestyle=dashed,dash=3pt 2pt]{*-*}(1.5,.5)(2.1,.5)
\psline[linestyle=dashed,dash=3pt 2pt]{*-*}(.1,.5)(.7,.5)
\psline[linestyle=dashed,dash=3pt 2pt]{-*}(.7,.5)(1.3,.3)
\psline[linestyle=dashed,dash=3pt 2pt]{-*}(.7,.5)(1.3,.7)
\end{pspicture}$$
The {\em degree} of such a diagram is 
half the number of all its vertices. 

Let $\CA_n^{\QQ}(M)$ denote the rational vector space generated by the degree $n$ diagrams
on $M$, quotiented out by the following relations AS and STU:
$$ {\rm AS:}  \begin{pspicture}[.2](0,-.2)(1,1)
\psline[linestyle=dashed,dash=3pt 2pt]{*-}(.5,.5)(.5,0)
\psline[linestyle=dashed,dash=3pt 2pt]{-}(.1,.9)(.5,.5)
\psline[linestyle=dashed,dash=3pt 2pt]{-}(.9,.9)(.5,.5)
\end{pspicture}
+
\begin{pspicture}[.2](0,-.2)(1,1)
\psline[linestyle=dashed,dash=3pt 2pt]{*-}(.5,.5)(.5,0)
\pscurve[linestyle=dashed,dash=3pt 2pt]{-}(.1,.9)(.7,.7)(.5,.5)
\pscurve[linestyle=dashed,dash=3pt 2pt]{-}(.9,.9)(.3,.7)(.5,.5)
\end{pspicture}
=0 \;\;\mbox{and} \;\;{\rm STU:} \begin{pspicture}[.2](0,-.2)(1,1)
\psline[linestyle=dashed,dash=3pt 2pt]{*-*}(.5,.5)(.5,.1)
\psline[linestyle=dashed,dash=3pt 2pt]{-}(.3,.9)(.5,.5)
\psline[linestyle=dashed,dash=3pt 2pt]{-}(.7,.9)(.5,.5)
\psline{->}(.1,.1)(1,.1)
\end{pspicture}
=
\begin{pspicture}[.2](0,-.2)(1,1)
\psline[linestyle=dashed,dash=3pt 2pt]{-*}(.3,.9)(.3,.1)
\psline[linestyle=dashed,dash=3pt 2pt]{-*}(.7,.9)(.7,.1)
\psline{->}(.1,.1)(1,.1)
\end{pspicture}
-
\begin{pspicture}[.2](0,-.2)(1,1)
\psline[linestyle=dashed,dash=3pt 2pt]{-*}(.3,.9)(.7,.1)
\psline[linestyle=dashed,dash=3pt 2pt]{-*}(.7,.9)(.3,.1)
\psline{->}(.1,.1)(1,.1)
\end{pspicture} 
$$
Each of these relations relate diagrams which can be represented by immersions that are identical outside the part of them represented in the pictures. For example, 
AS identifies the sum of two diagrams which only differ by the orientation
at one vertex to zero.

Set $\CA_n(M)=\CA_n^\RR(M) = \CA_n^{\QQ}(M) \otimes_{\QQ} \RR$ and let $$\CA(M) =\prod_{n \in \NN}\CA_n(M)$$ denote the product of the $\CA_n(M )$ as a topological vector space.
$\CA_0(M)$ is equal to $\RR$ generated by the empty diagram. The degree $n$
part of an element $\alpha=(\alpha_n)_{n\in \NN}$ of $\CA(M)$ will be denoted by $\alpha_n$.

\section{The Chern--Simons series}
\label{seccsi}

Let $M$ be an oriented one-manifold and let 
$$L : M \longrightarrow \RR^3$$ denote a $C^{\infty}$ embedding from $M$ to $\RR^3$. Let $\Gamma$ be a Jacobi diagram on $M$. Let $U=U(\Gamma)$ denote the set of univalent vertices of $\Gamma$, and let $T=T(\Gamma)$ denote the set of trivalent vertices of $\Gamma$. A {\em configuration\/} of $\Gamma$ is an embedding 
$$c:U \cup T \hookrightarrow \RR^3$$
whose restriction $c_{|U}$ to $U$ may be written as $L \circ j$ for some injection 
$$j:U \hookrightarrow M$$
in the given isotopy class $[i]$ of embeddings of $U$ into the interior of $M$. Denote the set of these configurations by
$C(L;\Gamma)$, $$C(L;\Gamma)=\left\{c:U \cup T \hookrightarrow \RR\sp{3} \; ; 
\exists j \in [i], c_{|U}=L \circ j\right\}.$$ 
In $C(L;\Gamma)$, the univalent vertices move along $L(M)$ while the trivalent vertices move in the ambient space, and $C(L;\Gamma)$ is naturally an open submanifold of 
$M^U \times (\RR^3)^T$.

Denote the set of (dashed) edges of $\Gamma$ by $E=E(\Gamma)$, and fix an orientation for these edges. Define the map
$\Psi:C(L;\Gamma) {\longrightarrow} \left(S\sp{2}\right)\sp{E}$
whose projection to the $S^2$ factor indexed by an edge from a vertex $v_1$ to a vertex $v_2$ is the direction of $\overrightarrow{c(v_1)c(v_2)}$.
This map $\Psi$ is again a map between two orientable manifolds that have the same dimension, namely the number of dashed half-edges of $\Gamma$, 
and we can write the {\em configuration space integral:\/}
$$I(L;\Gamma)=\int_{C(L;\Gamma)}\Psi\sp{\ast}(\Lambda\sp{E}\omega).$$
Bott and Taubes have proved that this integral is convergent \cite{bt}. 
Thus, 
this integral is well-defined up to sign. In fact, the orientation of the trivalent vertices of $\Gamma$ provides $I(L;\Gamma)$ with a well-defined sign.
Indeed, since $S^2$ is equipped with its standard orientation, it is enough to orient $C(L;\Gamma) \subset M^U \times (\RR^3)^T$ in order to define this sign. This will be done by providing the set of the natural coordinates of $M^U \times (\RR^3)^T$ with some order up to an even permutation. This set is in one-to-one correspondence with 
the set of (dashed) half-edges of $\Gamma$, and the vertex-orientation of the trivalent
vertices provides a natural preferred such one-to-one correspondence up to some (even!) cyclic permutations of three half-edges meeting at a trivalent vertex.
Fix an order on $E$, then the set of half-edges becomes ordered by
(origin of the first edge, endpoint of the first edge, origin of the second edge,
\dots, endpoint of the last edge), and this order orients $C(L;\Gamma)$.
The property of this sign is that the product $I(L;\Gamma)[\Gamma] \in {\cal A}(M)$ depends neither on our various choices
nor on the vertex orientation of $\Gamma$.

Now, the {\em perturbative series expansion of the Chern--Simons theory for one-manifold embeddings in $\RR^3$ \/} is the following sum running over all the Jacobi diagrams $\Gamma$ without
vertex orientation\footnote{This sum runs over equivalence classes of Jacobi diagrams, where two diagrams are equivalent if and only if they coincide except possibly for their vertex orientation.}:
$$Z_{\rm CS}(L)=
\sum_{\Gamma}{ \frac{I(L;\Gamma)}{\sharp \mbox{Aut}\Gamma}[\Gamma]} \;\; \in {\cal A}(M)$$
where $\sharp \mbox{Aut}\Gamma$ is the number of automorphisms of $\Gamma$ as 
a uni-trivalent graph with a given isotopy class of injections of $U$ into $M$, but
without vertex-orientation for the trivalent vertices.

Let $\theta$ denote the Jacobi diagram $$\theta=\tata$$ on $S^1$.
When $L$ is a knot $K$, the degree one part of $Z_{CS}(K)$ is $\frac{I(K;\theta)}{2}[\theta]$ and therefore $Z_{\rm CS}$ is not invariant under isotopy. However, the evaluation $Z^0_{\rm CS}$
at representatives of knots with null Gauss integral is an isotopy invariant that is a universal Vassiliev invariant of knots. (All the finite type knot invariants in the Vassiliev sense (see \cite{bn}) factor through it.) This is the content of
the following theorem, due independently to Altschuler and Freidel~\cite{af}, and to D.~Thurston~\cite{th}, after the work of many people including Guadagnini, Martellini and Mintchev~\cite{gua}, Bar-Natan~\cite{bn'}, Axelrod and Singer~\cite{as1,as2}, Kontsevich~\cite{Kon},  Bott and Taubes~\cite{bt}\dots

\begin{thm}[Altschuler--Freidel, D.~Thurston, 1995]
If $L=K_1\cup \dots \cup K_k$ is a link, then $Z_{CS}(L)$ only depends on the isotopy class of $L$ and on the Gauss integrals $I(K_i; \theta)$ of its components. In particular, the evaluation $$Z\sp{0}_{CS}(L) \in \prod_{n \in \NN}\CA_n(\sqcup_{i=1}^kS^1_i)$$ at representatives of $L$ whose components have zero Gauss integrals is an isotopy invariant of $L$. Furthermore, $Z\sp{0}_{CS}$ is a universal Vassiliev invariant of links.
\end{thm}

Recall that the normalized Kontsevich integral is also a universal Vassiliev knot invariant that is valued in the same target $\CA$. (See \cite{bn,les}.)
Thus, the still open natural question raised by Kontsevich in \cite{Kon} is:

{\em Is the Kontsevich integral of a zero-framed representative of a knot $K$ equal to $Z\sp{0}_{CS}(K)$?\/}

On one hand, this Chern--Simons series has a beautiful, very natural and completely symmetric definition. Furthermore, Dylan Thurston directly\footnote{The usual proof of the rationality of the Kontsevich integral of zero-framed links is indirect. It relies on the existence of a rational associator \cite{lm}.} proved that the Chern--Simons series is {\em rational\/}\footnote{For a link $L$, the degree $n$ part $Z\sp{0}_{CSn}(L)$ of  $Z\sp{0}_{CS}(L)$ belongs to $\CA^{\QQ}_n(\sqcup_{i=1}^kS^1_i)$.} because it behaves like a series of degrees of maps between closed manifolds \cite{th,Po}. 

On the other hand, the Kontsevich integral fits in with the framework of quantum link invariants and it can be defined in this setting \cite{lm,Ka}. Therefore, it is explicitly known how to recover quantum link invariants from the Kontsevich integral \cite{lm}. Furthermore, the computation of the Kontsevich integral for links can be reduced to the computation of small link pieces called {\em elementary q-tangles.\/} 
In \cite{Po}, Sylvain Poirier proved that the same can be done for the Chern--Simons series. Let us begin a review of his results.

\section{The Poirier extension of $Z_{CS}$ to tangles}
\label{secfunc}
 
A {\em planar configuration\/} is an embedding of a finite set $X$ into the  plane $\RR^2$.

In the {\em ambient space\/} $\RR^3=\{(x,y,z)\}$, the {\em horizontal plane\/} is the plane $(z=0)$, whereas the {\em blackboard plane\/} is the plane $(y=0)$. The $z$-coordinate
of a point $(x,y,z) \in \RR^3$ is called its {\em vertical projection\/}.

A {\em tangle\/} is the intersection of the image of a link representative transverse to $\RR^2 \times \{\beta,\tau\}$ with
a horizontal slice $\RR^2 \times [\beta,\tau]$. In particular, it is an embedded cobordism between two planar configurations. 

For $\lambda \in ]0,1]$, let $h_{\lambda}$ be the homeomorphism of $\RR^3$ that shrinks the horizontal plane with respect to the formula
$$h_{\lambda}(x_1,x_2,x_3)=(\lambda x_1,\lambda x_2,x_3).$$
Let $L(M)$ be a tangle. Set
$$Z_{P}(L(M))=\lim_{\lambda \rightarrow 0}Z_{CS}(h_{\lambda} \circ L) \in \CA(M).$$
Sylvain Poirier proved that this limit exists~\cite{Po}. 

Let $K$ be an almost horizontal knot embedding,
rotated by 90 degrees around a horizontal axis, (so that it is almost contained in some vertical plane) such that $I(K;\theta)=0$. Then $I(h_{\lambda} \circ K;\theta)=0$, for any $\lambda>0$, and the Poirier limit $Z_P$ is the limit of a constant map. Therefore, $Z_P$ is equal to the Chern--Simons series we started with for these representatives.

In general, we can see that the limit of $I(h_{\lambda} \circ K;\theta)$ depends on the differential degree of $\Psi$ near the equator of $S^2$. Assume that the height function (the third coordinate) of $K$ is a Morse function, (its second derivative does not vanish when the first one does). Then all the horizontal tangent vectors correspond to extrema of the height function. Identify the horizontal 
plane to $\CC$ so that the unit horizontal vector corresponding to an extremum $e$
is of the form $\exp(i\theta_e)$. When $\lambda$
approaches $0$, all the non-horizontal tangent vectors approach the poles, and the image of our annulus boundary (of Section~\ref{secgauss}) becomes a family of straight meridians intersecting the equator at directions $\exp(i\theta_e)$, $\exp(i(\theta_e+\pi))$ corresponding to extrema $e$. These meridians cut our sphere like an orange, and the differential degree
becomes constant on the boundaries of orange quarters and makes integral jumps at meridians. Thus, it can be seen that
$$\lim_{\lambda \rightarrow 0}I(h_{\lambda} \circ K;\theta) \equiv \frac{1}{\pi}\left(\sum_{\mbox{\footnotesize e minimum}} \theta_e - \sum_{\mbox{\footnotesize e maximum}} \theta_e\right) \mbox{mod}\; \ZZ$$ 
In particular, this limit is an integer when the horizontal vectors are in the blackboard plane, and does not vary under isotopies that keep the directions of the horizontal vectors fixed. Therefore, we define a {\em framed\/} tangle as a tangle whose horizontal vectors are contained in the blackboard plane. Two tangles are said to be {\em isotopic \/} if they can be obtained from one another by a composition of
\begin{itemize}
\item a possible {\em rescaling of the height parameter,\/} that is a composition by $1_{\RR^2}
\times h$ where $h$ is an increasing diffeomorphism from $[\beta,\tau]$ to another interval of $\RR$, and
\item an isotopy that is framed at any time, and that preserves the bottom and the top planar configurations up to homotheties with positive ratio and translations. 
\end{itemize}
The homotheties with positive ratio will be called {\em dilations.\/}

In \cite{Po}, Sylvain Poirier proved that his limit $Z_P$ is a well-defined functorial isotopy invariant of framed tangles. He also interpreted his limit for braids (that are paths of planar configurations) as the Chen holonomy of a connection on a trivial bundle over the space of planar configurations. 
He proved that, unlike the holonomy of the complex Knizhnik-Zamolodchikov connection
that provides a similar definition of the Kontsevich integral for braids \cite{bn,les}, the holonomy of his real connection converges without any regularisation for paths reaching limit planar configurations. Thus, he was able
to extend his invariant to framed cobordisms between limit\footnote{These limits are precisely defined in \cite[Section 10.1]{Po}
as the boundary points of a suitable compactification of the space of planar configurations, where the planar configurations are considered up to translations and dilations.
In this compactification, points are allowed to collide with each other. But the compactification provides us with magnifying glasses allowing us to see any restricted configuration at the scales of the collisions.}
configurations.

The Le and Murakami {\em framed q-tangles\/} of \cite{lm} are framed cobordisms between particular
limit linear configurations. In~\cite{lm}, Le and Murakami extended the Kontsevich integral to a monoidal functor of framed $q$-tangles that behaves
"nicely" under the operations of deleting a component, doubling a component that runs from bottom to top, reversing orientations, and under the orthogonal symmetries
whose fix point sets are either the vertical axis, or the blackboard plane or the horizontal plane.
The framed q-tangles, these operations and the expected "nice" behaviour are precisely defined in \cite{Les}. We call a functor that satisfies all these properties
of the Kontsevich--Le--Murakami functor a {\em good\footnote{Again, the complete definition can be found in \cite{Les}.}\/} functor.

The general properties that Sylvain Poirier proved for his limit $Z_P$ in \cite{Po} imply the following theorem.

\begin{thm} The restriction of $Z_P$ to framed q-tangles is  a good functor.
\end{thm}

Furthermore, Poirier expressed the variation of $Z_P$ under a {\em framing change\/}, where $\smaxo$ becomes $\smaxno$, in terms of 
the following constant called the {\em anomaly.\/}

\section{The anomaly}

Let us define the anomaly. Let $v \in S^2$. Let $D_v$ denote the linear map
$$\begin{array}{llll}D_v: &\RR &\longrightarrow& \RR^3\\
& 1&\mapsto&v\end{array}$$
Let $\Gamma$ be a Jacobi diagram on $\RR$.
Define $C(D_v;\Gamma)$ and $\Psi$ as in Section~\ref{seccsi}.
Let $\hat{C}(D_v;\Gamma)$ be the quotient of $C(D_v;\Gamma)$ by the translations parallel to $D_v$ and by the dilations. Then $\Psi$ factors through
$\hat{C}(D_v;\Gamma)$ that has two dimensions less. Now, allow $v$ to run through $S^2$ and define $\hat{C}(\Gamma)$ as the total space
of the fibration over $S^2$ where the fiber over $v$ is $\hat{C}(D_v;\Gamma)$. The map $\Psi$ becomes a map between two smooth oriented
manifolds of the same dimension. Indeed, $\hat{C}(\Gamma)$ carries a natural smooth structure and can be oriented as follows. Orient $C(D_v;\Gamma)$ as before, orient $\hat{C}(D_v;\Gamma)$ so that $C(D_v;\Gamma)$ is locally homeomorphic to the oriented product
(translation vector $(0,0,z)$ of the oriented line, ratio of homothety $\lambda \in ]0,\infty[$) $\times \hat{C}(D_v;\Gamma)$ and orient $\hat{C}(\Gamma)$ with the $(\mbox{base} (=S^2) \oplus \mbox{fiber})$ convention\footnote{This can be summarized by saying that the $S^2$-coordinates replace $(z,\lambda)$.}.  Then we can again define 
$$I(\Gamma)=\int_{\hat{C}(\Gamma)}\Psi\sp{\ast}(\Lambda\sp{E}\omega).$$
Now, the {\em anomaly\/}\footnote{Our convention differs from the Poirier convention of \cite{Po} where $\alpha$ denotes half this element.} is the following sum running over all connected Jacobi diagrams $\Gamma$ on the oriented lines (again without vertex-orientation):
$$\alpha=\sum{ \frac{I(\Gamma)}{\sharp \mbox{Aut}\Gamma}[\Gamma]} \;\; \in {\cal A}(\RR).$$
Its degree one part is 
$$\alpha_1= \left[ \begin{pspicture}[0.4](0,0)(.4,.8)
\psline{->}(0.05,0.05)(0.05,0.75)
\pscurve[linestyle=dashed,dash=2pt 1pt]{*-*}(0.05,.15)(.3,.35)(.05,.5)
\end{pspicture} \right]. $$
The central symmetry of $\RR^3$ acts on $\hat{C}(\Gamma)$ (by composition of the configurations). Studying the different
orientation changes induced by this action shows that, for any $\Gamma$, $$I(\Gamma)[\Gamma]=(-1)^{1+E(\Gamma)+T(\Gamma)} I(\Gamma)[\Gamma]=(-1)^{n+1} I(\Gamma)[\Gamma]$$
Therefore, for any integer $n$,
$$\alpha_{2n}=0.$$
In \cite{lm}, Le and Murakami proved that any 
good functor that varies like the Kontsevich integral $Z_K$ under a framing change
coincides with $Z_K$ for framed links. 
This theorem allowed Poirier to derive the following corollary from his work.

\begin{thm}[Poirier \cite{Po}]
If the anomaly $\alpha$ vanishes in degree greater than one, then the Kontsevich
integral of zero framed links is equal to
their Chern--Simons series $Z^0_{CS}$.
\end{thm}
Say that an element $\beta=(\beta_n)_{n \in \NN}$ in $\CA(S^1)$ is a {\em two-leg element\/} if, for any $n \in \NN$, $\beta_n$ is a combination of diagrams with two univalent 
vertices.
Let $\beta$ be a two-leg element.
Forgetting $S^1$ from $\beta$ gives rise to a unique series $\beta^s$ of diagrams
with two distinguished univalent vertices $v_1$ and $v_2$, such that $\beta^s$ is symmetric with respect to the exchange of $v_1$ and $v_2$. The series $\beta^s$ is well-defined thanks to the diagrammatic Bar-Natan version \cite{bn} of the Poincar\'e-Birkhoff-Witt theorem.\\
A {\em chord diagram\/} on a one-manifold $M$ is a Jacobi diagram without trivalent vertices. Its (dashed connected) components are just {\em chords.\/} The degree $n$ chord diagrams generate $\CA_n(M)$.
If $\Gamma$ is a chord diagram, define $\Phi(\beta)([\Gamma])$ by replacing each chord
by $\beta^s$. Then $\Phi(\beta)$ is a well-defined morphism of topological vector spaces from $\CA(M)$ to $\CA(M)$ for any one-manifold $M$, and  $\Phi(\beta)$ is an isomorphism as soon as $\beta_1 \neq 0$. See \cite{Les}.

In \cite{Les}, I refined the Le and Murakami uniqueness theorem, by characterizing the possible variations of good functors under framing changes. This allowed me to refine the above Poirier theorem as follows.

\begin{thm}\label{thmLes}{\rm\cite{Les}}\qua
The anomaly $\alpha$ is an odd\footnote{$\alpha_{2n}=0$.} two-leg element of $\CA(S^1)$.
For any framed link $L$, the Poirier limit integral $Z_P(L)$ is equal to
$\Phi(\alpha)(Z_K(L))$.
Thus, for any zero-framed link $L$, $Z^0_{CS}(L)$ is equal to
$\Phi(\alpha)(Z_K(L))$.
\end{thm}

The above result yields an algebraic constraint on the anomaly. Since, conversely, any functor of the form $\Phi(\alpha)(Z_K(L))$, for a real odd two-leg element $\alpha$, is a good functor, nothing more can be obtained from algebra. We are now going to try to compute the low degree terms of the anomaly in a geometric way.

\section{Some geometric properties of the anomaly}
\label{secgano}

Poirier showed that the anomaly can be defined from the logarithm of the holonomy of his connexion for the two-strand braid \begin{pspicture}[0.4](0,0)(.4,.6)
\psecurve{-}(.2,-.15)(.35,0)(.2,.15)(.05,.3)(.2,.45)
\psecurve[border=1.5pt](.2,-.15)(.05,0)(.2,.15)(.35,.3)(.2,.45)(.05,.6)(.2,.75)
\psecurve[border=1.5pt](.2,.15)(.05,.3)(.2,.45)(.35,.6)(.2,.75)
\end{pspicture}. 
This is restated in the following proposition.

Let $\Gamma$ be a Jacobi diagram on $\{1,2\} \times \RR$.
For $(x,y) \in \RR^2$, let $$\begin{array}{llll}f_{(x,y)}:
& \{1,2\} \times \RR &\longrightarrow & \RR^3\\
&(1,t) & \mapsto &(0,0,t)\\
&(2,t) & \mapsto &(x,y,-t) \end{array}$$
For $\theta \in [0,2\pi[$, let $C_{\theta}(\Gamma)$ denote the quotient of 
$C(f_{(\cos(\theta),\sin(\theta))};\Gamma)$ by the translations by some $(0,0,z)$, with $z \in \RR$.
And let $C(\Gamma)=\cup_{\theta \in [0,2\pi[}C_{\theta}(\Gamma)$. Orient $C(\Gamma)$ by letting the 
$\theta$-coordinate 
$\theta \in [0,2\pi[$ replace\footnote{This means that if the quotient $C_{\theta}(\Gamma)$ is oriented so 
that $C(f_{(\cos(\theta),\sin(\theta))};\Gamma)$ is oriented by the (fiber $\oplus$ base) convention where the base is $C_{\theta}(\Gamma)$ 
and the fiber is the oriented vertical translation factor $\RR$ of $\RR^3$, then $C(\Gamma)$ is oriented with the (base$=[0,2\pi[$ $\oplus$ fiber $=C_{\theta}(\Gamma)$ ) convention.} the translation parameter 
$z \in \RR$.

Define the {\em two-strand anomaly\/} $$\tilde{\alpha}=
\sum_{\Gamma \;\mbox{\footnotesize connected Jacobi diagram on } \{1,2\} \times \RR}\frac{1}{\sharp \mbox{Aut}(\Gamma)}
\int_{C(\Gamma)}\Psi^{\ast}(\Lambda^{E(\Gamma)}\omega)[\Gamma].$$
Let $i: \CA(\trid) \longrightarrow \CA(\RR)$ be the linear continuous map induced 
by the inclusion from $\{1,2\} \times \RR =\trid$ to $ \begin{pspicture}[.2](0,0)(.5,.6)
\pscurve{->}(.05,0)(0.15,0.6)(.25,0)\end{pspicture}
=
\begin{pspicture}[.2](0,0)(.5,.6)
\psline(.05,0)(.05,.4)
\pscurve[linestyle=dashed,dash=2pt 1pt](.05,.4)(0.15,0.6)(.25,.4)
\psline{->}(.25,.4)(.25,0) \end{pspicture}$.
The map $i$ sends a Jacobi diagram $\Gamma$ to the diagram with the same dashed graph
equipped with the same orientations at trivalent vertices
where the embedding of univalent vertices is composed by the above inclusion.

\begin{proposition}[Poirier \cite{Po}]
\label{proptwos}
$$\alpha=-i(\tilde{\alpha}).$$
\end{proposition}

This definition of the anomaly is easier to handle with. Instead of 
searching for configurations where the univalent vertices are on a common 
unknown line, we look for configurations where the univalent vertices have 
the same horizontal coordinates.

\begin{proposition}
\label{mainprop} Let $n$ be an integer greater than 2.
Let $\Gamma$ be a degree $n$ connected Jacobi diagram on $\{1,2\} \times \RR$.
Set $$I(\Gamma)=\int_{C(\Gamma)}\Psi^{\ast}(\Lambda^{E(\Gamma)}\omega).$$
If $\Gamma$ has less than three vertices on some strand, then
$I(\Gamma)=0.$

In particular, the two-strand anomaly $\tilde{\alpha}_n$ in degree $n$ 
is a combination of connected diagrams 
with at least 3 vertices on each strand.
\end{proposition}
\bp
Let $\Gamma$ be a diagram as in the above statement.
Let $U_i$ denote the set of its univalent vertices which are on $\{i\} \times \RR$, for $i \in \{1,2\}$.
We shall prove that if $U_2$ contains less than 3 elements, then  
$I(\Gamma)=0.$

This will be sufficient to conclude by symmetry. Indeed,
if $\Gamma$ is such that $U_1$ contains less than 3 elements, then the diagram $s(\Gamma)$ obtained from $\Gamma$ by exchanging the two lines is such that $U_2$ contains less than 3 elements; and, by considering the action of the central symmetry of $\RR^3$,
by composition of the configurations (up to translations), we can see that $I(\Gamma)=\pm I(s(\Gamma))$.

If $U_2$ is empty, then all the $(2\sharp E(\Gamma)-1)$-manifolds $C_{\theta}(\Gamma)$ are the same.
Thus $\Psi(C(\Gamma))$ is the image under $\Psi$ of one of them, and its volume in 
$\left(S^2\right)^{E(\Gamma)}$ is zero.

In order to study the two remaining cases where $U_2$ contains one or two elements, we replace
$C(\Gamma)$ by a slightly larger smooth configuration space $\hat{C}(\Gamma)$ in which
$C(\Gamma)$ is dense, where $\Psi$ extends so that we shall have
$$I(\Gamma)=\int_{\hat{C}(\Gamma)}\Psi^{\ast}(\Lambda^{E(\Gamma)}\omega).$$
Fix one vertex $u_0$ on $\{1\} \times \RR$, and let $C_0(f_{(x,y)};\Gamma)$
denote the subset of $C(f_{(x,y)};\Gamma)$ made of the configurations that map $u_0$ to the origin of $\RR^3$. Orient it as the quotient of $C(f_{(x,y)};\Gamma)$ by the vertical translations with the ((fiber $(=\RR)$) $\oplus$ base)
convention.
Let $\hat{C}(\Gamma)$
denote the $(2\sharp E(\Gamma))$-dimensional quotient of 
$$P(\Gamma)=\cup_{(x,y) \in \RR^2}C_0(f_{(x,y)};\Gamma)$$
by the dilations with a ratio $\lambda \in ]0,\infty[$.
The map $\Psi$ is well-defined on this space which contains the additional configurations corresponding
to $(0,0)$ that constitute a codimension 2 subspace of $\hat{C}(\Gamma)$ that will therefore not contribute to the integral.

The orientation of $\hat{C}(\Gamma)$ is obtained as follows. Orient $P(\Gamma)$ with the convention $(\mbox{base}= \RR^2) \oplus (\mbox{fiber})$. Then the orientation of $\hat{C}(\Gamma)$ is defined
by the $(\mbox{fiber}= ]0,\infty[) \oplus (\mbox{base})$ convention. 

Now, let us get rid of the case where $U_2=\{u\}$, by defining a free smooth action of
$]0,\infty[$ on $\hat{C}(\Gamma)$ that does not change the image of a configuration under
$\Psi$, so that $\Psi$ will again factor through a map from a $(2 \sharp E-1)$-dimensional manifold to $(S^2)^E$ and define a zero integral.
Let $t$ denote the other end of the edge of $u$. Let $\mu \in ]0,\infty[$. Let $c \in P(\Gamma)$. Define 
$\mu.c(x)=x$ if $x \neq u$ and $\mu.c(u)=c(t)+\mu(c(u)-c(t))$. 
This action is compatible with the action of the dilations (and different since the degree of $\Gamma$ is greater that $1$), $\Psi(\mu.c)=\Psi(c)$, and we are finished with this case.

Let us study the remaining case, $U_2=\{u_1,u_2\}$. Denote the trivalent vertex connected by an edge to $u_1$ by $t_1$,
and denote the trivalent vertex connected by an edge to $u_2$ by $t_2$. The vertices $t_1$ and $t_2$ may coincide.
We shall use the symmetry\footnote{This symmetry resembles 
a Bott and Taubes symmetry that is used to prove the isotopy invariance of $Z^0_{CS}$.}
$\sigma$ that maps a configuration $c \in P(\Gamma)$ to the configuration 
$\sigma(c)$ defined by:\\
$$\begin{array}{lll}\sigma(c)(x)&=c(x) & \mbox{if}\; x \notin U_2\\
\sigma(c)(u_1)&=c(t_1)+c(t_2)-c(u_2)&\\
\sigma(c)(u_2)&=c(t_1)+c(t_2)-c(u_1)&\end{array}$$ 
as shown in the picture below.

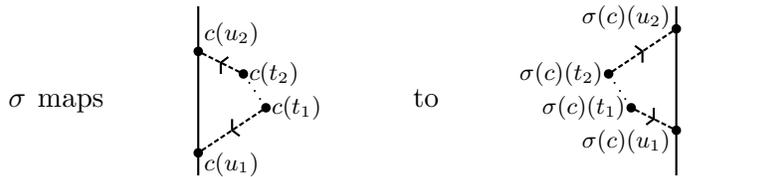
\begin{center}

$\sigma$ maps $\begin{pspicture}[.4](-1.1,0)(2.7,2.3)
\psset{yunit=1.5cm}
\psline(0,0)(0,1.5)
\rput[l](.95,.6){\footnotesize $c(t_1)$}
\rput[l](.65,.9){\footnotesize $c(t_2)$}
\rput[l](.05,1.25){\footnotesize $c(u_2)$}
\rput[l](.05,.1){\footnotesize $c(u_1)$}
\psline(.4,1.05)(.3,1)(.3,.9)
\psline(.55,.35)(.45,.4)(.45,.5)
\psline[linestyle=dashed,dash=2pt 1pt]{*-*}(.9,.6)(0,.2)
\psline[linestyle=dashed,dash=2pt 1pt]{*-*}(.6,.9)(0,1.1)
\psline[linestyle=dotted](.6,.9)(.9,.6)
\end{pspicture}$ to
$\begin{pspicture}[.4](-1.5,0)(2.7,2.3)
\psset{yunit=1.5cm}
\psline(1.5,0)(1.5,1.5)
\rput[r](.85,.6){\footnotesize $\sigma(c)(t_1)$}
\rput[r](.55,.9){\footnotesize $\sigma(c)(t_2)$}
\rput[r](1.45,1.4){\footnotesize $\sigma(c)(u_2)$}
\rput[r](1.45,.3){\footnotesize $\sigma(c)(u_1)$}
\psline(1.1,.45)(1.2,.5)(1.2,.6)
\psline(.95,1.15)(1.05,1.1)(1.05,1)
\psline[linestyle=dashed,dash=2pt 1pt]{*-*}(.9,.6)(1.5,.4)
\psline[linestyle=dashed,dash=2pt 1pt]{*-*}(.6,.9)(1.5,1.3)
\psline[linestyle=dotted](.6,.9)(.9,.6)
\end{pspicture}$.

\end{center}

This symmetry factors through the dilations and reverses the orientation
of $\hat{C}(\Gamma)$.
Furthermore, $\Psi(\sigma(c))$ is obtained from $\Psi(c)$ by reversing 
the two unit vectors corresponding to the edges containing $u_1$ and $u_2$ and by exchanging them,
that is, by a composition by a diffeomorphism of $\left(S^2\right)^{E(\Gamma)}$ that preserves $\Lambda^{E(\Gamma)}\omega$.
This proves that $I(\Gamma)$ vanishes and this finishes the proof of the proposition.
\eop

As a corollary, $\alpha_n$ is a combination of Jacobi diagrams with 6 legs when $n \geq 2$.
Unfortunately, unlike the previous two-leg condition, this condition is not very restrictive.
Nevertheless, this provides another proof\footnote{Note that with this two-strand definition all configuration space integrals vanish 
whereas with the one-strand original definition, the different integrals cancel each 
other witout being zero individually. See~\cite[Section~7]{Po}.} of
the Poirier and Yang result that $\alpha_3=0$. There is no connected degree 3 diagram with at least 6
univalent vertices.   Recall that for any integer $n$, $\alpha_{2n}=0$. 
Therefore, the next interesting degree is 5. The connected degree 5 diagrams with at 
least 6 univalent vertices are necessarily trees and their dashed parts have one of the two forms:

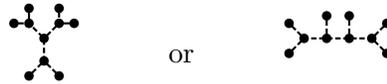
\begin{center}
$\begin{pspicture}[.2](-1.1,0)(1.9,.9)
\psline[linestyle=dashed,dash=2pt 1pt]{*-*}(0,.7)(.2,.7)(.4,.5)(.4,.2)(.2,0)
\psline[linestyle=dashed,dash=2pt 1pt]{*-*}(.4,.2)(.6,0)
\psline[linestyle=dashed,dash=2pt 1pt]{*-*}(.2,.7)(.2,.9)
\psline[linestyle=dashed,dash=2pt 1pt]{*-*}(.4,.5)(.6,.7)(.6,.9)
\psline[linestyle=dashed,dash=2pt 1pt]{*-*}(.6,.7)(.8,.7)
\end{pspicture}$ 
or
$\begin{pspicture}[.2](-1.1,0)(2.4,.9)
\psline[linestyle=dashed,dash=2pt 1pt]{*-*}(0,.3)(.2,.5)(1.1,.5)(1.3,.3)
\psline[linestyle=dashed,dash=2pt 1pt]{*-*}(.5,.5)(.5,.8)
\psline[linestyle=dashed,dash=2pt 1pt]{*-*}(.8,.5)(.8,.8)
\psline[linestyle=dashed,dash=2pt 1pt]{*-*}(0,.7)(.2,.5)
\psline[linestyle=dashed,dash=2pt 1pt]{*-*}(1.1,.5)(1.3,.7)
\end{pspicture}$ 
\end{center}

Now, observe that if $\Gamma$ is a Jacobi diagram on $\{1,2\} \times \RR$, 
with two univalent vertices that are connected to the same trivalent vertex, and that lie on the same vertical
line, like in
$\begin{pspicture}[.2](-.1,0)(.6,1)
\psline[linestyle=dashed,dash=2pt 1pt]{*-}(0,.2)(.3,.5)(.5,.5)
\psline[linestyle=dashed,dash=2pt 1pt]{*-*}(0,.8)(.3,.5)
\psline(0,0)(0,1)
\end{pspicture}$, then $I(\Gamma)$ vanishes.
Indeed the two corresponding edges and the vertical vector are coplanar.
Therefore, the image of $\Psi$ must lie inside a codimension one subspace of $\left(S^2\right)^{E(\Gamma)}$.
This additional remark determines the distribution of the univalent vertices of the above graphs on the two 
vertical lines. Then 
Sylvain Poirier \cite{Pop} computed $\tilde{\alpha}_5$
with the help of Maple, and he found that $\alpha_5=0$ thanks to AS and STU.
As a corollary, all coefficients of the HOMFLY polynomial properly normalized that are Vassiliev invariants of degree less than seven can be explicitly written as 
combinations of the configuration space integrals of Section~\ref{seccsi}. 

Thus, the following Bar-Natan theorem generalizes to any {\em canonical\/} \footnote{Here, {\em canonical\/} can be understood as explicitly recovered from the Kontsevich integral like all the quantum invariants.} Vassiliev invariant of degree less than 7. 

\begin{thm}[Bar-Natan \cite{bn'}, 1990]
Let $\Delta$ denote the symmetrized Alexander polynomial. For any knot $K$, 
$$\frac{\Delta^{\prime \prime}(K)(1)}2=-\frac13I(K;\begin{pspicture}[0.4](-.7,-.6)(.7,.6)
\psarc{->}(0,0){.6}{180}{20}
\psarc{-}(0,0){.6}{20}{180}
\SpecialCoor
\psline[linestyle=dashed,dash=2pt 1pt]{*-*}(.6;-30)(0,0)(.6;-150)
\psline[linestyle=dashed,dash=2pt 1pt]{*-*}(0,0)(.6;90)
\end{pspicture})+\frac14I(K;\begin{pspicture}[0.4](-.7,-.6)(.7,.6)
\psarc{->}(0,0){.6}{180}{20}
\psarc{-}(0,0){.6}{20}{180}
\SpecialCoor
\psline[linestyle=dashed,dash=2pt 1pt]{*-*}(.6;-45)(.6;135)
\psline[linestyle=dashed,dash=2pt 1pt]{*-*}(.6;-135)(.6;45)
\end{pspicture})+\frac1{24}.$$
\end{thm}

This particular coefficient, that is the degree 2 invariant which can be extracted from the Chern--Simons series, has been further studied in \cite{pv}.


\section{Questions}

The first question here is of course:

(1)\qua Prove or disprove the physicist conjecture:\\
\begin{center}
$\alpha_i=0$ for any $i>1$.
\end{center}

Let us briefly summarize what remains to be done here. The space $\CA(\RR)=\CA(S^1)$ splits as a natural direct sum where the space of two-leg diagrams is a direct summand, thanks to the diagrammatic Bar-Natan version \cite{bn} of the Poincar\'e-Birkhoff-Witt theorem. Call the natural projection on this summand  the {\em two-leg projection.\/}
According to Theorem~\ref{thmLes} and Proposition~\ref{proptwos}, it remains to compute the two-leg projection of $i(\tilde{\alpha}_{2n+1})$ for $n \geq 3$, where
Proposition~\ref{mainprop} and Lemma~1.9 in \cite{Po} allow us to neglect numerous diagrams in this computation.

After the articles of Axelrod, Singer~\cite{as1,as2}, Bott and Cattaneo~\cite{bc1,bc2,cat},
Greg Kuperberg and Dylan Thurston have constructed a universal finite type
invariant for homology spheres as a series of configuration space integrals
similar to $Z^0_{CS}$, in \cite{KT}. Their construction yields two natural questions:

(2)\qua Find a surgery formula for the Kuperberg-Thurston invariant in terms of the above
Chern--Simons series.

(3)\qua Compare the  Kuperberg--Thurston invariant to the LMO invariant \cite{lmo}.

\Addresses\recd
\end{document}

%% file: gtmonout.tex
\def\ifplaintex{\expandafter\ifx\csname documentclass\endcsname\relax}

\ifplaintex 
\else
\fi

\def\gtm{{\mathsurround=0pt\it $\cal G\mskip-2mu$eometry \&\ 
$\cal T\!\!$opology $\cal M\mskip-1mu$onographs}}    

\def\gtp{{\mathsurround=0pt\it $\cal G\mskip-2mu$eometry \&\ 
$\cal T\!\!$opology $\cal P\!$ublications}}  

\def\recd{{\small Received:\qua\receiveddate\ifx\reviseddate\relax
\else\qquad Revised:\qua\reviseddate\fi\par}} 

\def\volumenumber#1{\def\thevolumenumber{#1}}
\def\volumeyear#1{\def\thevolumeyear{#1}}
\def\volumename#1{\def\thevolumename{#1}}
\def\papernumber#1{\def\thepapernumber{#1}}
\def\pagenumbers#1#2{\def\startpage{#1}\def\finishpage{#2}}
\def\published#1{\def\publishdate{#1}}
\def\received#1{\def\receiveddate{#1}}
\def\revised#1{\def\reviseddate{#1}}
\def\accepted#1{\def\accepteddate{#1}}

\def\asciiaddress#1{\def\theasciiaddress{#1}}

\long\def\asciiabstract#1{\long\def\theasciiabstract{#1}}
\def\asciikeywords#1{\def\theasciikeywords{#1}}

\let\\\par
\let\thevolumenumber\relax\let\thepapernumber\relax
\let\thevolumeyear\relax\let\startpage\relax
\let\finishpage\relax\let\publishdate\relax\let\receiveddate\relax
\let\reviseddate\relax\let\accepteddate\relax\let\theasciititle\relax
\let\theasciiauthors\relax\let\theasciiaddress\relax
\let\theasciiabstract\relax\let\theasciikeywords\relax

\let\theerratum\relax\let\theasciiemail\relax
\let\theshortauthors\relax\let\theshorttitle\relax

\def\startpage{1}\def\finishpage{15}\def\thepapernumber{77}

\volumenumber{2}
\volumename{Proceedings of the Kirbyfest}
\volumeyear{1999}

\long\def\maketitlep{   

\count0=\startpage

\gtm\nl        
{\small Volume \thevolumenumber: \thevolumename\nl 
\ifx\theerratum\relax\else Erratum \erratumnumber\nl\fi
Pages \startpage--\finishpage\nl}

\vglue 0.1truein   

{\parskip=0pt\leftskip 0pt plus 1fil\def\\{\par\smallskip}{\ifplaintex\large
\else\Large\fi\bf\thetitle}\par\medskip}   
\vglue 0.05truein 

{\parskip=0pt\leftskip 0pt plus 1fil\def\\{\par}{\sc\theauthors}
\par\medskip}%
 
\vglue 0.03truein 

{\small\leftskip 25pt\rightskip 25pt{\bf Abstract}\stdspace\theabstract

{\bf AMS Classification}\stdspace\theprimaryclass
\ifx\thesecondaryclass\relax\else; \thesecondaryclass\fi\par
{\bf Keywords}\stdspace \thekeywords\par}\vglue 7pt

}   

\font\phead=cmsl9 scaled 950
\font=cmsl9 scaled 1050
\font\pnum=cmbx10 scaled 913
\font\lnum=cmbx10 
\font\pfoot=cmsl9 scaled 950
\font=cmsl9 scaled 1050
\ifplaintex
\headline{\vbox to 0pt{\vskip -4.5mm\line{\small\phead\ifnum
\count0=\startpage ISSN 1464-8997 (on line)
1464-8989 (printed) \hfill {\pnum\folio}\else\ifodd\count0\def\\{ }%
\ifx\theshorttitle\relax\thetitle\else\theshorttitle\fi\hfill{\pnum\folio}
\else\def\\{ and }{\pnum\folio}\hfill\ifx\theshortauthors\relax\theauthors
\else\theshortauthors\fi\fi\fi}\vss}}
\footline{\vbox to 0pt{\vglue 0mm\line{\small\pfoot\ifnum\count0=\startpage
Published \publishdate:\qua\copyright\ \gtp\hfill\else
\gtm, Volume \thevolumenumber\ (\thevolumeyear)\hfill\fi}\vss
}}
\else
\makeatletter
\def\@oddhead{{\small\lhead\ifnum\count0=\startpage ISSN 1464-8997 (on line)
1464-8989 (printed) \hfill {\lnum\number\count0}\else\ifodd\count0
\def\\{ }\ifx\theshorttitle\relax \thetitle \else\theshorttitle\fi\hfill
{\lnum\number\count0}\else\def\\{ and }{\lnum\number\count0}
\hfill\ifx\theshortauthors\relax 
\theauthors\else\theshortauthors\fi\fi\fi}}\def\@evenhead{@oddhead}
\def\@oddfoot{\small\lfoot\ifnum\count0=\startpage Published \publishdate:\qua\copyright\ \gtp\hfill\else
\gtm, Volume \thevolumenumber\ (\thevolumeyear)\hfill\fi}
\def\@evenfoot{@oddfoot}
\makeatother
\fi

\let\maketitlepage\maketitlep
\let\makeshorttitle\maketitlepage
\let\maketitle\maketitlepage

\newwrite\gtoutfile
\long\gdef\makeheadfile{  
{\def\\{, }\def\s{ }
\immediate\openout\gtoutfile head.xxx
\immediate\write\gtoutfile{To: math@arxiv.org}
\immediate\write\gtoutfile{Subject: put OR rep NNNNN:ppppp}
\immediate\write\gtoutfile{--text follows this line--}
\immediate\write\gtoutfile{Proxy-for: \ifx\theasciiauthors\relax
\theauthors\else\theasciiauthors\fi\s<\ifx\theasciiemail\relax\theemail\else\theasciiemail\fi>}
\immediate\write\gtoutfile{\noexpand\\}
\immediate\write\gtoutfile{Authors: \ifx\theasciiauthors\relax
\theauthors\else\theasciiauthors\fi}
{\def\\{ }\immediate\write\gtoutfile{Title: \ifx\theasciititle\relax
\thetitle\else\theasciititle\fi}}
\immediate\write\gtoutfile{Subj-class: GT or SG, GR etc}
\immediate\write\gtoutfile{MSC-class: \theprimaryclass\ifx\thesecondaryclass\relax\else, \thesecondaryclass\fi}
\immediate\write\gtoutfile{Journal-ref: Geom. Topol. Monogr. \thevolumenumber\s
(\thevolumeyear) \startpage-\finishpage}
\immediate\write\gtoutfile{Comments: Published by Geometry and Topology Monographs at}
\immediate\write\gtoutfile{\s\s\s  http://www.maths.warwick.ac.uk/gt/GTMon\thevolumenumber/paper\thepapernumber.abs.html}
\immediate\write\gtoutfile{\noexpand\\}
\immediate\write\gtoutfile{}
\ifx\theasciiabstract\relax
\immediate\write\gtoutfile{\theabstract}\else
\immediate\write\gtoutfile{\theasciiabstract}\fi
\immediate\write\gtoutfile{}
\immediate\write\gtoutfile{\noexpand\\}
\immediate\write\gtoutfile{}
\immediate\closeout\gtoutfile}}  

\def\maketitlepage{\maketitlep\makeheadfile}
\let\makeshorttitle\maketitlepage
\let\maketitle\maketitlepage